\documentclass[12pt,a4paper,draft]{article}

\usepackage[reqno]{amsmath}
\usepackage{amssymb,amsthm,upref,amscd}
\usepackage[T1]{fontenc}
\usepackage{times}

\newtheorem{theorem}{Theorem}[section]
\newtheorem{corollary}[theorem]{Corollary}
\newtheorem{definition}[theorem]{Definition}

\newtheorem{lemma}[theorem]{Lemma}
\newtheorem{proposition}[theorem]{Proposition}
\newtheorem{remark}[theorem]{Remark}

\theoremstyle{definition} \theoremstyle{remark}
\numberwithin{equation}{section}

\begin{document}

\title{The existence of  solutions for a Schr\"{o}dinger equation with jumping nonlinearities crossing the essential spectrum }
\author{Chong Li$^{\mathrm{a,b,c}}$\thanks{%
The first author is supported by NSFC(11871066), NSFY(Y911251) and
NSFC(11471319). Email:
lichong@amss.ac.cn.}\ and \ \ Xinyu Li$^{\mathrm{d}}$\thanks{%
Corresponding author. Email:
lixinyu@amss.ac.cn.}\  \\
\\
{\small $^{\mathrm{a}}$Institute of Mathematics, Academy of Mathematics and
Systems Science, }\\
{\small Chinese Academy of Sciences, Beijing 100190, China }\\
{\small $^{\mathrm{b}}$Hua Loo-Keng Key Laboratory of Mathematics, Chinese
Academy of Sciences }\\
{\small $^{\mathrm{c}}$Center for Excellence in Mathematical Sciences,
Chinese Academy of Sciences (CEMS, CAS)}\\
{\small $^{\mathrm{d}}$Department of Mathematics, Tsinghua University, Beijing, 100084, PR China
 }}
\date{}
\maketitle

\begin{abstract}
In this paper, we establish the existence of one solution  for a Schr\"{o}dinger equation with jumping nonlinearities:
$-\Delta u+V(x)u=f(x,u)$, $x\in \mathbb {R}^N$, and $u(x)\to 0$, $|x|\to +\infty$,
where $V$ is a potential function on which we make hypotheses, and in particular allow $V$ which is unbounded below, and $f(x,u)=au^-+bu^++g(x,u)$. No restriction on $b$ is required, which implies that $f(x,s)s^{-1}$ may interfere with the essential spectrum of $ -\Delta+V$ for $s\to +\infty$.
Using the truncation method and the Morse theory, we can compute  critical  groups of the  corresponding functional at zero and infinity, then  prove the existence of   one negative  solution.

\bigskip



\noindent\textbf{Data availability statement:} Our manuscript has no associate data.
\end{abstract}


\medskip

\section{Introduction}
This paper is mainly concerned with nonlinear Schr\"{o}dinger equations of the form
\begin{align}\label{2}
\begin{cases}
-\Delta u+V(x)u=f(x,u), & x\in \mathbb {R}^N,\\
u(x)\to 0, & as\quad |x|\to +\infty,
\end{cases}
\end{align}
arising from study of standing wave solutions of time-dependent nonlinear Schr\"{o}dinger equations.
The corresponding energetic functional of (\ref{2}) is

\begin{equation*}
J(u)=\frac12\int_{\mathbb{R}^N}|\nabla u|^2+V(x)u^2 dx-\int_{\mathbb{R}^N} F(x,u),\nonumber
\end{equation*}
$F(x,u)=\int ^u_0 f(x,s)ds$,
where
$f(x,u)$ is a Carath\'{e}odory function on $\mathbb {R}^N\times\mathbb{R}$ such that
\begin{align}
\frac {f(x,t)}t\to \begin{cases}
a & as  \quad t \to -\infty,\\
b & as \quad t \to +\infty,
\end{cases}
\end{align}
$(a,b)\notin \sum(-\Delta +V)$, and $\sum(-\Delta +V)$ is called Fu\v{c}\'{i}k spectrum, defined as the set of all $(a,b)\in \mathbb{R}^2$ such that
\begin{align}\label{3}
\begin{cases}
-\Delta u+V(x)u=au^-+bu^+ & x\in \mathbb {R}^N,\\
u(x)\to 0  & as \quad |x|\to +\infty,
\end{cases}
\end{align}
has a nontrivial solution $u$ in the form domain of $-\Delta +V$, where $u^+=\max \{u,0\}$, $u^-=\min \{u,0\}$. We denote by $\sigma(-\Delta+V)$ the spectrum of $-\Delta+V$, and by $\sigma_{ess}(-\Delta+V)$ its essential spectrum.
The form domain of $-\Delta +V$ is the Hilbert space
\[
H^1_m(\mathbb{R}^N):=\{ u\in H^1(\mathbb{R}^N): \int_{\mathbb{R}^N}|\nabla u|^2+(V+m) u^2< +\infty \},
\]
endowed with the scalar product
\[
\langle u,v\rangle_m := \int_{\mathbb{R}^N} \nabla u\nabla v+(V+m)uv  \quad\text{for} \quad u,v\in H^1_m(\mathbb{R}^N),
\]
and induced norm
\[
\|u\|_m:= \left( \int_{\mathbb{R}^N}| \nabla u|^2+(V+m)u^2\mathrm {d}x\right) ^{\frac 12},
\]
where $m> -\inf\sigma(-\Delta+V)$ is arbitrary but fixed.

In order to show the existence of a negative solution of (\ref{2}), we assume on the potential $V$:\\
$(V_1)$ Let $V$ be real-valued,  $V=V^++V^-$,   $V^+\in L^r_{loc}(\mathbb{R}^N)$ and $V^-$ is a Kato-Rellich potential, i.e.,  $V^-(x)\in L^p(\mathbb{R}^N)+L^{\infty}(\mathbb{R}^N)$,  $r, p=2$  if $N\leq 3$, $r, p> 2$ if $N=4$, $r, p>\frac N2$ if $N\geq 5$;\\
$(V_2)$ $\lambda _1=\inf \sigma(-\Delta+V)<\lambda_2<\cdots<\lambda_{l-1}<\lambda_l\leq\sigma_0:=\inf\sigma_{ess}(-\Delta+V), l\geq 2$;\\
And also on the nonlinear term $f$:\\
($f_1$) Set $f(x,s)=as^-+bs^++g(x,s)$, $g(x,0)=0$, $\lim\limits_{|s|\to\infty}\frac{g(x,s)}s=0$, uniformly with respect to $x\in \mathbb{R}^N$. Moreover, there exists $\beta >0$  such that for  any  $s_1$, $s_2\in \mathbb{R}$ and $ x\in\mathbb{R}^N$, s.t.
\begin{equation}\label{33}
|g(x,s_1)-g(x,s_2)|\leq \beta |s_1-s_2|,
\end{equation}
$\beta <\sigma _0-a$;\\
($f_2$) Set $f(x,s)=a_0s^-+b_0s^++\tilde g(x,s),\lim\limits_{s\to 0_-}\frac{f(x,s)}s=a_0$, $\lim\limits_{s\to 0_+}\frac{f(x,s)}s=b_0$ uniformly on $x\in \mathbb{R}^N$.

The main features of problem (\ref{2}) is that the nonlinearity is asymptotically linear and the potential $V$ is allowed to be unbounded below.
Our main result is the following theorem.
\begin{theorem} \label{them1} Assume that $(V_1)$, $(V_2)$, $(f_1)$, and $(f_2)$ hold. If   $a>\lambda_1 >a_0$,
then (\ref{2}) admits one negative solution.
\end{theorem}
\begin{remark}
Note that Theorem \ref{them1} continues to hold  when   $V=V_1+V_2$, where $V$ is bounded below, $V_1\geq 0$, $ V_1\in L^2_{loc}(\mathbb{R}^N)$ and  $V_2$ is a $-\Delta$-bounded multiplication operator with relative bound less than one.
\end{remark}
\begin{remark}
Under ($V_1$) and $\inf\limits_{x\in\mathbb{R}^N} V(x)> -\infty$,
the results  of the  existence of four solutions of Schr\"{o}dinger  equation and the constructions of minimal and maximal curves of Fu\v{c}\'{i}k spectrum in \cite{12}, are still true, since the  assumption ($V_1$) implies unique continuation properties(see \cite{14,22,23}for related results).
\end{remark}

Motivations for the study  of problems with jumping nonlinearities come from problems of  Mechanics and Engineering, as floating beam equations, models of suspension bridges, etc... (see \cite{27,28,29} and the references therein).

The Fu\v{c}\'{i}k spectrum was originally introduced on the smooth bounded domain $\Omega \in\mathbb{R}^N$ by Fu\v{c}\'{i}k \cite{8} and Dancer \cite{7} in the 1970s, which plays an important role in the study of some elliptic problems with linear growth. Existence and multiplicity of solutions for problems of this type are strictly related to the position of pair $(a,b)$ with respect to the Fu\v{c}\'{i}k spectrum. In fact, these problems may lack compactness in the sense that the well known Palais-Smale compactness condition fails if the pair $(a,b)$ belongs to the Fu\v{c}\'{i}k spectrum $\Sigma$; moreover, the topological properties of the sublevels of the corresponding energy functional depend on the position of $(a,b)$ with respect to $\Sigma$.

In \cite{3} and \cite{12}, the authors dealt with the Fu\v{c}\'{i}k spectrum for Schr\"{o}dinger operator on $\mathbb{R}^N$.
For results  concerning the study of multiplicity of solutions  with  asymptotically linear term or jumping nonlinearities, we mention, for example \cite{2,12,13,14,15} (and references therein). In \cite{14}, the authors obtained a positive solution, a negative solution and a sign-changing solution by  using the invariant sets of descending flow. See also \cite{6,25,26,31,32} for further studies on problems in bounded set with jumping nonlinearities or problems where the nonlinear terms interfere with the spectrum of the linear operator
and
Schr\"{o}dinger type equations with an asymptotically linear term.

To apply variational methods, it is important to know whether Palais-Smale sequences are compact. The compactness of such sequences strongly depends on the interplay of the nonlinearity $f$ with the spectrum of $-\Delta +V$.
There are very few works on asymptotically linear problems on unbounded domains when the essential spectrum $\sigma_{ess}(-\Delta+V)$ exists.

In \cite{3,23}, compactness of Palais-Smale sequences was observed  for equation (\ref{2}) with potentials $ V_{\lambda}=a+\lambda g$ which have a potential well whose steepness is controlled by a real parameter $\lambda$. As $\lambda\to \infty$, the potential $ V_{\lambda}$ has a  well, and $\lim_{\lambda\to \infty}\inf{\sigma_{ess}}(-\Delta+V_{\lambda})=\infty$.

In \cite{12}, the authors obtained four solutions of Schr\"{o}dinger equation problem with $V$  a Kato-Rellich potential,  by verifying the   (PS) condition,  constructing of minimal and maximal curves of Fu\v{c}\'{i}k spectrum in $Q_l$, building up a minimax principle on the form domain of $-\Delta+V$,   and  finally proving a weak maximum principle for $\mathbb{R}^N$, which  serves as a tool to get two critical points in positive and negative cone, respectively.

The case $(a^*,b^*]\cap \sigma_{ess}(-\Delta + V)\ne \varnothing$
tends to be more intricate, where $a^*=\inf\limits_{s\in[0,+\infty)}\frac{f(x,s)}s$ and  $b^*=\sup\limits_{s\in[0,+\infty)}\frac{f(x,s)}s$. Usually, additional monotonicity assumptions are  imposed on $V$ or $f$. In \cite{19}, the authors obtained one positive radial solution with $V\equiv a$   constant, where $f(x,s)$ is radially symmetric on $x$,  and $f(r,s)s^{-1}$ is a non-decreasing function of $s$ on $[0,\infty)$. In \cite{11}, the existence of a positive
solution was proved, where $V$ is a constant, $f(x,s)s^{-1}\to a\in(0,\infty]$ as $s\to +\infty$, $f(\cdot,s)$ is periodic in $x_i$, $1\leq i\leq N$ and $F(x,s)s^{-2}$ is non-decreasing function of $s\geq 0$.

In \cite{24,11} the proof of the boundedness of Palais-Smale sequences relies on the work of Lions on the concentration compactness principle.
\cite{24} established the existence of a positive solution for an asymptotically linear elliptic problem on $\mathbb {R}^N$ autonomous at infinity with $V$ a continuous potential.
It is however more delicate, in contrast to just periodic in \cite{11}, because of the spectral structure of $-\Delta+V$ and the fact that the equation  is not enjoying a translation invariance.

In contrast to these works, in \cite{12} the authors do not use the monotonicity, but instead  take advantage of the position of $(a,b)$ which is below the $\inf \sigma_{ess}(-\Delta+V)$. Key to our approach to  prove  the Palais-Smale condition  is the use of results on problems established by \cite{12}.

In \cite{14} $V$ and $f$ are allowed to depend on $x$ without any assumption concerning periodicity or the existence of a limit for $|x|\to \infty$, and the nonlinearity is allowed to interact locally with the essential spectrum of $-\Delta +V$. The compactness is proved for Palais-Smale sequences or Cerami sequences at any energy level.

In this work, we extend the argument by \cite{12} and allow $V$ is unbounded below.
Our objective is to find nontrivial solutions of (\ref{2}), without imposing either a monotonicity  assumption on $f$ or a bounded condition on $V$.
We observe that no restriction on the global behavior of $V$ is required other than $(V_1)$. In particular, $V$ is not required to be bounded or to belong to a Kato class. On the other hand, the hypotheses on the nonlinearity $f$ are fairly mild, in the sense that
we assume $f(x,u)=au^-+bu^++g(x,u)$  without  restricting  $b$, which implies that $f(x,s)s^{-1}$ may interfere with the essential spectrum of $ -\Delta+V$ for $s\in (0,+\infty)$.

We will define a suitable modification of equation (\ref{2}) of which the corresponding  functional $J$ satisfies the Palais-Smale  condition and for which we will find a critical point via Morse theory. This critical point will eventually be shown to be  a negative solution of the original equation, taking advantage of the weak maximum principle for $\mathbb{R}^N$ in \cite{12}.

The problem of the  self-adjointness of $-\Delta+V$  has a long history, for which we refer the reader to a  book by Reed and Simon \cite{17}.
The general theory of relatively bounded perturbations,  for example
the Kato-Rellich theorem,  will permit us to consider the Schr\"{o}dinger operators  of the form  $-\Delta+V$
with a wide class of potentials $V$, where
$V=V^++V^-$  is real-valued,  $V^+\in L^r_{loc}(\mathbb{R}^N)$ and $V^-$ is a Kato-Rellich potential, with $r=2$ if $N\leq 3$,  $r>\frac N2$ if $N\geq 4$.

Close analysis of such potentials   finally leads us to  capture  the properties of the infimum of the spectrum, which turns out to be crucial  in the computation of  our main results in the critical  groups of the  corresponding functional at infinity.
Indeed,
we will  show that the eigenvector associated to an eigenvalue strictly above the infimum of the spectrum of $-\Delta+V$ with such potentials must be sign changing, using the results that the discrete spectrum of $-\Delta+V$ is determined by the behavior of $-\Delta+V$ on bounded subsets of $ \mathbb {R}^N$ and the essential spectrum of $-\Delta+V$ is determined by the behavior of $-\Delta+V$  ( in particular, $V$) in a neighborhood of infinity.     As a corollary, when the  infimum of the spectrum is a point spectrum, the corresponding eigenvector must be non-negative (or non-positive).

Indeed,
that the bottom of the spectrum of $-\Delta+V$ is a simple eigenvalue  $\lambda_1$
was already assumed, and the only remaining problem was to show that its corresponding  eigenvector $\varphi_1$ of $\lambda_1$ guarantees a given sign.

In the study of quantum mechanical energy operators, particular attention is paid to the question of whether the bottom of the spectrum is an eigenvalue and whether that eigenvalue is simple.
To derive the property of the eigenvector $\varphi_1$,
we benefit from the results on the non-degenerate ground state for  Schr\"{o}dinger operators.
Suppose that $V\in L^1_{loc}(\mathbb{R}^N)$  with $V\geq 0$ and let $H=-\Delta +V$ be defined as the sum of quadratic forms.  If  $\inf \sigma(H)$ is an eigenvalue, then it is a simple eigenvalue and the corresponding eigenfunction is strictly positive (see proposition \ref{prop31}, and \cite{18,30} for related results).

Since the potential $V$ may be unbounded from below, we can not use the result above directly.
Our approach proposed here consists of  constructing a sequence of potentials $\{V_n\}$ bounded below such that  the bottom of the spectrum of $-\Delta+V_n$ is simple.

By assumption $(V_1)$, we set $V^-=\tilde V^*_1+\tilde V^*_2$, where $ \tilde V^*_1\in L^p(\mathbb{R}^N)$, $\tilde V^*_2\in L^{\infty}(\mathbb{R}^N)$, $p=2$ if $N\leq 3$,  $p>\frac N2$ if $N\geq 4$.
We introduce a potential sequence $\{V_n\}$, where $V_n=V^++\tilde V^*_{1,n}+\tilde V^*_2$,  $\tilde V^*_{1,n}=
\tilde V^*_1$ if  $ \tilde V^*_1\geq -n$, otherwise $\tilde V^*_{1,n}=-n$.
It is obvious that $ V_n$ is bounded below, then  the result  for the ground state can apply to $-\Delta+ V_n$, if we could prove that $\lambda_1^{(n)}<\sigma_0^{(n)}$ for $n$ large enough, where $\lambda_1^{(n)}$ and $\sigma_0^{(n)}$ denote the simple eigenvalue and bottom of the essential spectrum of $-\Delta+V_n$  respectively.

In the spectrum theory, Persson (see \cite{10,16}) discovered a beautiful geometric description for the bottom of the essential spectrum of a semi-bounded Schr\"{o}dinger operator. Indeed,
\begin{equation}\label{1471}
\inf \sigma_{ess}(H)=\sup _{\kappa\subset\mathbb{R}^N}\inf _{\phi\neq 0,\phi \in C^{\infty}_0(\mathbb{R}^N\setminus \kappa)}\frac {\langle\phi ,H\phi\rangle}{\|\phi\|^2_{L^2}},
\end{equation}
where the supremum is over all compact subset $ \kappa \subset \mathbb{R}^N$.
On the other hand, in order to distinguish the lowest eigenvalue with the essential spectrum, we need the following equation
\begin{equation}
\inf_{u\in H^1_m,\|u\|_{L^2}=1}\int_{\mathbb{R}^N}|\nabla u|^2+V u^2,
\end{equation}
which expresses the lowest eigenvalue of a self-adjoint operator as the infimum of a quadratic form.

The fact that $V^-$ belongs to Kato-Rellich class  and  Persson's formula guarantee that $\lambda_1^{(n)}\to \lambda_1$ and $\sigma_0^{(n)}\to \sigma_0$ as $n \to \infty$,  which  gives us the desired result.
Then denote by $\varphi_1^{(n)}$ the corresponding eigenfunction of $\lambda_1^{(n)}$, with $\varphi_1^{(n)}>0$.  The fact that $\varphi_1^{(n)}$ is an exact critical point often provides additional information on the sequence $\{\varphi_1^n\}$, which helps to show its convergence.

The idea of constructing Palais-Smale sequences that possess some extra properties  that might  be helpful, is an old topic.
Let $I(u)=\frac 12 \int_{\mathbb{R}^N}|\nabla u|^2  +V(x) u^2 -\frac 12\lambda_1 \int_{\mathbb{R}^N}u^2 dx$, and
$I_n(u)=\frac 12 \int_{\mathbb{R}^N}|\nabla u|^2  +V_n(x) u^2 -\frac 12\lambda_1^{(n)}\int_{\mathbb{R}^N}u^2 dx$.
Instead of starting from a sequence of approximate critical points of $I$ (as in the case of a standard Palais-Smale sequence), we now start from a sequence of exact critical points of nearby functionals.
We have  $I_n^{\prime}(\varphi_1^{(n)})= 0$ for all $n\in  \mathbb{N}$ and $I^{\prime}(\varphi_1^{(n)})\to 0$, as $n \to \infty$. We will show $\varphi_1^{(n)}\to \varphi_1$.  Then we conclude that $\varphi_1\geq 0$. This is done in Lemma \ref{LM1}.

At this point, we will reduce  the computation of  critical groups   of  $J$ at zero and infinity to  the  computation  of critical groups of corresponding  homogeneous functional at zero, inspired by the paper \cite{12}.

The  organization of this paper is as follows: In section \ref{sec2}, we give some Preliminaries result and functional setting. In Section \ref{sec3}, we study the properties of the infimum of the spectrum and compute  critical  groups of the  corresponding functional at zero and infinity. Finally, Section \ref{sec4} is devoted to the proof of Theorem \ref{them1}.

$\bf{Notation.}$
Throughout the article,
we use the following $L^p(\mathbb {R}^N)$ and $H^1(\mathbb {R}^N)$-norms:
\[ \|u\|^p_{L^p(\mathbb{R}^N)}=\int_{\mathbb{R}^N}|u|^p dx,  \forall p\in (1,\infty),
\]
\[ \|u\|_{H^1(\mathbb{R}^N)}^2=\int_{\mathbb{R}^N}|\nabla u|^2+ u^2 dx,
\]
where there is no risk of confusion, we will denote $\|u\|^p_{L^p(\mathbb{R}^N)}$ simply by $\|u\|_{L^p}$ and $\|u\|_{H^1(\mathbb{R}^N)}^2$ by $\|u\|_{H^1}$.
Also if we take a subsequence of a sequence $\{u_n\}$ we shall denote it again $\{u_n\}$.
\section{Preliminaries}\label{sec2}
\subsection{Self-adjointness of $-\Delta+V$}
The problem of the essential self-adjointness of $-\Delta+V$  has a long history, for which we refer the reader to a  book by Reed and Simon \cite{17}.

\begin{definition}(see page 136, \cite{10})
let $V(x)\in L^p(\mathbb{R}^N)+L^{\infty}(\mathbb{R}^N)$, and be real, with $p=2$ if $N\leq 3$, $p> 2$, if $N=4$ and $p>\frac N2$, if $N\geq 5$, then $V$ is called Kato-Rellich (K-R) potential.
\end{definition}

Enormous background knowledge is involved in
self-adjointness. What we are concerned about is, if $A$ is self-adjoint and $B$ is hermitian, under which conditions is $A+B$ self-adjoint?
To show the self-adjointness of $-\Delta+V$, we need the following theorems.
\begin{definition}(see page 162, \cite{17})
Let $A$ and $B$ be densely defined linear operators on a Hilbert space $H$. Suppose that:

$(i)$ $D(B)\supset D(A)$,

$(ii)$ For some $a$ and $b$ in $\mathbb{R}$ and all $\varphi\in D(A)$,
\begin{equation}\label{AB}
\|B\varphi\|\leq a\|A\varphi\|+b\|\varphi\|.
\end{equation}
Then $B$ is said to be $A$-bounded. The infimum of such $a$ is called the relative bound of $B$ with respect to $A$.

\end{definition}

\begin{proposition}
(the Kato-Rellich Theorem, see Theorem \uppercase\expandafter{\romannumeral 10}.12, \cite{17}).
Suppose that $A$ is  self-adjoint, $B$ is symmetric, and  $B$ is $A$-bounded with relative bound $a<1$. Then $A+B$ is self-adjoint on $D(A)$ and essentially self-adjoint on any core of $A$. Further, if $A$ is bounded below by $M$, then $A+B$ is  bounded below by $M-\max\{\frac b{1-a},b+a|M|\}$ where $a$ and $b$ are given by (\ref{AB}).

\end{proposition}
The results that K-R potential is  $-\Delta$-bounded and the relative bound of K-R potential with respect to $-\Delta$ is $0$, are concluded by applying
 the Kato-Rellich theorem  to perturbations of the free particle hamiltonian
\begin{equation*}
dom H_0=H^2(\mathbb {R}^N),     H_0=-\Delta.
\end{equation*}

\begin{proposition}(see Theorem IX.28, \cite{17})
Let $\varphi\in L^2(\mathbb{R}^n)$ be in $D(H_0)$. Then

$(a)$ If $n\leq 3$, $ \varphi$ is a bounded continuous function and for any $a>0$, there is a $b$, independent of $\varphi$, so that
\begin{equation*}
\|\varphi\|_{L^{\infty}}\leq a \|H_0\varphi\|+ b\|\varphi\|,
\end{equation*}

$(b)$ If $n\geq 4$, and $2\leq q<\frac{2n}{n-4}$, then $\varphi\in L^q(\mathbb{R}^n)$ and for any $a>0$, there is a $b$ (depending only on $q$, $ n$ and $a$) so that

\begin{equation*}
\|\varphi\|_{L^q}\leq a \|H_0\varphi\|+ b\|\varphi\|.
\end{equation*}
\end{proposition}

For more general potentials $V$,
to show the self-adjointness of $-\Delta+V$,  we need the following theorems.

\begin{proposition}(see Theorem X.28, \cite{17})
Let $V\in L^2_{loc}(\mathbb{R}^N)$ with $V\geq 0$ pointwise. Then $-\Delta+V$ is essentially self-adjoint on $C_0^{\infty}(\mathbb{R}^N)$.
\end{proposition}

\begin{proposition}(the Davies-Faris theorem, see Theorem X.31, \cite{17}) Let $H_0$ be a positive self-adjoint operator on $L^2(X,d\mu)$ such that $e^{-tH_0}$ is positively preserving. Let $V$ be a multiplication operator with $V\geq 0$. Suppose that $H=H_0+V$ is essentially self-adjoint on $D(H_0)\cap D(V)$. Let $W$ be an $H_0$-bounded multiplication operator. Then $W$ is $H$-bounded. In fact, if
\begin{equation}\label{H0}
\|W\psi\|\leq a \|(H_0+b)\psi\| \quad\forall \psi\in D(H_0),
\end{equation}
then
\begin{equation}\label{H}
\|W\psi\|\leq a \|(H+b)\psi\| \quad\forall \psi\in D(H).
\end{equation}

\end{proposition}

As a corollary of the Davies-Faris theorem we have
\begin{corollary}
Let $V=V_1+V_2$ where $V_1\geq 0$,$ V_1\in L^2_{loc}(\mathbb{R}^N)$ and where $V_2$ is a $-\Delta$-bounded multiplication operator with relative bound $a<1$. Then $-\Delta+V$ is essentially self-adjoint on $C_0^{\infty}(\mathbb{R}^N)$.
\end{corollary}
\begin{remark}
Let $V$ be a multiplication operator with $V^+\in L^2_{loc}(\mathbb{R}^N)$ and $V^-\in  L^p(\mathbb{R}^N)+L^{\infty}(\mathbb{R}^N)$,   $p=2$ if $N\leq 3$, $p> 2$ if $N=4$, $p>\frac N2$ if $N\geq 5$. Then $-\Delta+V$ is essentially self-adjoint on $C_0^{\infty}(\mathbb{R}^N)$. By the the Kato-Rellich theorem, $D(-\Delta+V)=D(-\Delta+V^+)$.
\end{remark}
\begin{remark}\label{remk5} Let $V$ be a multiplication operator with $V^+\in L^2_{loc}(\mathbb{R}^N)$ and $V^-\in  L^p(\mathbb{R}^N)+L^{\infty}(\mathbb{R}^N)$,  where $p=2$ if $N\leq 3$, $p> 2$, if $N=4$ and $p>\frac N2$, if $N\geq 5$. Since $-\Delta+V$ is essentially self-adjoint on $C_0^{\infty}(\mathbb{R}^N)$, and bounded below,  by choosing $m>0$ suitably large, clearly, $A_m:=-\Delta+V+m$ is a positive definite self-adjoint operator with $D(-\Delta+V+m)=D(-\Delta+V^++m):=H^2_m(\mathbb{R}^N)=\{u\in L^2(\mathbb{R}^N): (-\Delta+V^+)u\in  L^2(\mathbb{R}^N)\}=\{u\in L^2(\mathbb{R}^N): (-\Delta+V)u\in  L^2(\mathbb{R}^N)\}$, and the form domain  $H^1_m(\mathbb{R}^N)=\{ u\in H^1(\mathbb{R}^N): \int_{\mathbb{R}^N}|\nabla u|^2+(V+m) u^2< +\infty \}=\{ u\in H^1(\mathbb{R}^N): \int_{\mathbb{R}^N}|\nabla u|^2+(V^++m) u^2< +\infty \}$.
\end{remark}

\subsection{The equivalence of space norms}
Let $V$ be a real-valued measurable function with $V=V^++V^-$, where  $V^+\in L^r_{loc}(\mathbb{R}^N)$ and $V^-$ is a K-R potential, i.e. $V^-(x)\in L^p(\mathbb{R}^N)+L^{\infty}(\mathbb{R}^N)$,  $r,p=2$ if $N\leq 3$, $r,p> 2$ if $N=4$, $r,p>\frac N2$ if $N\geq 5$.

It might just as well assume $m+\lambda_1>0$ and $\inf\sigma(\hat A)>1$ with $\hat A=-\Delta+2V^-+m$ by choosing $m>1$ suitably large. Clearly, $A_m:=-\Delta+V+m$ is a positive definite self-adjoint operator on the Hilbert space $D(A_m)=H^2_m(\mathbb{R}^N)=\{u\in L^2(\mathbb{R}^N): (-\Delta+V)u\in  L^2(\mathbb{R}^N)\}$.

Denote by
\begin{equation*}
H^1_m(\mathbb{R}^N):=\{ u\in H^1(\mathbb{R}^N): \int_{\mathbb{R}^N}|\nabla u|^2+(V+m) u^2< +\infty \},
\end{equation*}
equipped with the scalar product
\begin{equation*}
\langle u,v\rangle_m:= \int_{\mathbb{R}^N} \nabla u\nabla v+ (V+m)uv,
\end{equation*}
\\and the induced norm
\begin{equation*}
\|u\|_m := \left( \int_{\mathbb{R}^N}| \nabla u|^2+(V+m)u^2\mathrm {d}x\right) ^{\frac 12}.
\end{equation*}
Obviously,
$H^1_m(\mathbb{R}^N)\subset H^1(\mathbb{R}^N)$ via

\begin{align*}
\frac 12\|u\|_{H^1}^2&= \frac12\int_{\mathbb{R}^N}| \nabla u|^2+u^2 dx\nonumber\\
           &\leq \frac 12\int_{\mathbb{R}^N}| \nabla u|^2+\frac 12 \int_{\mathbb{R}^N}| \nabla u|^2+V^+u^2+(2V^-+m)u^2 dx\nonumber\\
           &\leq \|u\|^2_m.
\end{align*}
Define
\begin{equation*}
\|u\|_{*,m}= \left( \int_{\mathbb{R}^N}| \nabla u|^2+(V^++m)u^2\mathrm {d}x\right) ^{\frac 12}.
\end{equation*}
It is obvious that $\|u\|_{*,m}\geq \|u\|_{H^1}$, since we choose $m>1$ suitably large.

Indeed, notice that
$V^-$ is a K-R potential, set $V^-=\tilde V^*_1+\tilde V^*_2$, $\tilde V^*_1\in L^p(\mathbb{R}^N)$, $\tilde V^*_2\in L^{\infty}(\mathbb{R}^N)$, and be real, where $p=2$ if $N\leq 3$, $p> 2$, if $N=4$ and $p>\frac N2$, if $N\geq 5$,
then for any $u\in H^1_m(\mathbb{R}^N)$,
\begin{equation*}
\int_{\mathbb{R}^N}|V^-u^2|\leq (C_0\|\tilde V^*_1\|_{L^p}+  \|\tilde V^*_2\|_{L^{\infty}})\|u\|_{H^1}^2.
\end{equation*}

Divide the proof into two cases:\\
(i) $N\geq 4$. Since $p>\frac N2$,$ \frac 1p +\frac 1q=1\Rightarrow q\in (1,\frac {N}{N-2})$. By H\"{o}lder inequality and Cauchy-Schwarz inequality,
\begin{align*}
\int_{\mathbb{R}^N}|V^-u^2|&\leq \|\tilde V^*_1\|_{L^p}\|u\|^2_{L^{2q}} + \|\tilde V^*_2\|_{L^{\infty}}\|u\|^2_{L^2}\nonumber\\
                        &\leq C_0\|\tilde V^*_1\|_{L^p}\|\nabla u\|^2_{L^2}+ \|\tilde V^*_2\|_{L^{\infty}}\|u\|^2_{L^2}\nonumber\\
                        &\leq(C_0\|\tilde V^*_1\|_{L^p}+  \|\tilde V^*_2\|_{L^{\infty}}+1)\|u\|_{H^1}^2.
\end{align*}
(ii) $N\leq 3$.  Take $p=2$. Similarly,
\begin{align*}
\int_{\mathbb{R}^N}|V^-u^2|&\leq \|\tilde V^*_1\|_{L^p}\|u\|^2_{L^{2q}} + \|\tilde V^*_2\|_{L^{\infty}}\|u\|^2_{L^2}\nonumber\\
                        &\leq C_0\|\tilde V^*_1\|_{L^2}\|\nabla u\|^2_{L^2}+ \|\tilde V^*_2\|_{L^{\infty}}\|u\|^2_{L^2}\nonumber\\
                        &\leq(C_0\|\tilde V^*_1\|_{L^2}+  \|\tilde V^*_2\|_{L^{\infty}}+1)\|u\|_{H^1}^2.
\end{align*}
Observe that
\begin{align*}
\frac 12\|u\|^2_{*,m}&=\frac12 \int_{\mathbb{R}^N}| \nabla u|^2+(V^++m)u^2 dx\nonumber\\
                     &\leq\frac12 \int_{\mathbb{R}^N}| \nabla u|^2+(V^++m)u^2 dx+\frac 12 \int_{\mathbb{R}^N}| \nabla u|^2+V^+u^2+(2V^-+2m) u^2 dx\nonumber\\
           &\leq \int_{\mathbb{R}^N}| \nabla u|^2+(V+m)u^2 dx =\|u\|_m^2.
\end{align*}
On the other hand,
\begin{align*}
\|u\|_m^2&= \int_{\mathbb{R}^N}| \nabla u|^2+(V+m)u^2 dx\nonumber\\
           &\leq \int_{\mathbb{R}^N}| \nabla u|^2+(V^++m)u^2+\int_{\mathbb{R}^N}|V^-|u^2 dx\nonumber\\
           &\leq \|u\|^2_{*,m}+(C_0\|\tilde V^*_1\|_{L^p}+  \|\tilde V^*_2\|_{L^{\infty}})\|u\|_{H^1}^2\nonumber\\
           &\leq (C_0\|\tilde V^*_1\|_{L^p}+  \|\tilde V^*_2\|_{L^{\infty}}+1)\|u\|^2_{*,m}.
\end{align*}
The assertion follows.

\section{The existence of one negative solution}\label{sec3}

\subsection{Truncation trick}\label{3.1}
Let $V$ be a real potential as above, then
there exists $m\in \mathbb{R}$, such that $A_m=-\Delta +(V+m)$ is a positive definite self-adjoint operator. In what follows, we assume

($f_1$) Set $f(x,s)=as^-+bs^++g(x,s)$, $g(x,0)=0$, $\lim\limits_{|s|\to\infty}\frac{g(x,s)}s=0$, uniformly with respect to $x\in \mathbb{R}^N$. Moreover, there exists $\beta >0$ such that for any  $s_1$, $s_2\in \mathbb{R}$, and $ x\in\mathbb{R}^N$ there holds
\[
|g(x,s_1)-g(x,s_2)|\leq \beta |s_1-s_2|
\]
and $\beta <\sigma _0-a$.

($f_2$) Set $f(x,s)=a_0s^-+b_0s^++\tilde g(x,s)$, $\lim\limits_{s\to 0_-}\frac{f(x,s)}s=a_0$, $\lim\limits_{s\to 0_+}\frac{f(x,s)}s=b_0$ uniformly on $x\in \mathbb{R}^N$.

Set
\begin{align}
f^*_m(x,t)= \begin{cases}
(a+m)t+g(x,t) & as  \quad t <0,\\
0 & as \quad t\geq 0.
\end{cases}
\end{align}
Since there exists $\beta >0$ such that for any  $s_1$, $s_2\in \mathbb{R}$, and $ x\in\mathbb{R}^N$ there holds
\[
|g(x,s_1)-g(x,s_2)|\leq \beta |s_1-s_2|\nonumber.
\]
Suppose $\beta<a+m$, then we have $f^*_m(x,u)\leq 0$.

If there exists a non-positive solution of
\begin{align}\label{fix}
\begin{cases}
-\Delta u+(V+m)u=f^*_m(x,u)=(a+m)u^-+g(x,u^-),& x\in \mathbb {R}^N,\\
u(x)\to 0,& as |x|\to +\infty,
\end{cases}
\end{align}
with the aid of the weak maximum principle established by \cite{12}, we get the solution is negative, so it is also a negative solution of (\ref{2}) of the form
\begin{equation}\label{fm}
\begin{cases}
-\Delta u+(V+m)u=(a+m)u^-+(b+m)u^++ g(x,u),& x\in \mathbb {R}^N,\\
u(x)\to 0,& as |x|\to +\infty.
\end{cases}
\end{equation}
\\
We associate with (\ref{fix}) the functional $J: H^1_m(\mathbb{R}^N) \to \mathbb{R}$  defined by

\[
  J(u)=\frac12\int_{\mathbb{R}^N}|\nabla u|^2+(V(x)+m)u^2 dx-\frac 12\int_{\mathbb{R}^N} (a+m)(u^-)^2-\int_{\mathbb{R}^N} G^*(x,u),
\]
where
$G^*(x,u)=\int ^{u^-}_0 g(x,s)ds$.
It is noted that the weak solutions of (\ref{fix}) are the critical points of $J$ in $H^1_m(\mathbb{R}^N)$.

On account of the hypothesis, $\sigma (A_m)\in (0,+\infty)$ alludes to the fact that $A_m^{-1}$ exists. Obviously, in view of $(f_1)$, $A_m^{-1}f^*_m(x,u)\in H^2_m(\mathbb{R}^N)$ for any $u\in H^1_m(\mathbb{R}^N)$.
For fixed $\varphi\in H^1_m(\mathbb{R}^N),$
\begin{align*}
\int_{\mathbb{R}^N}f^*_m(x,u)\varphi&=\int_{\mathbb{R}^N}A_mA_m^{-1}f^*_m(x,u)\varphi\nonumber\\
                                  &=\nabla A_m^{-1}f^*_m(x,u) \nabla \varphi+(V+m)A_m^{-1}f^*_m(x,u)\varphi.
\end{align*}
\\Therefore,
\begin{align*}
\langle J^{\prime}(u),\varphi\rangle_m&=\langle u-A_m^{-1}f^*_m(x,u),\varphi\rangle_m\nonumber\\
                                      &=\int_{\mathbb{R}^N} \nabla u\nabla\varphi +(V(x)+m)u\varphi -\int_{\mathbb{R}^N}f^*_m(x,u)\varphi,
\end{align*}
yielding $J^{\prime}(u)=u-A_m^{-1}f^*_m(x,u)$.

In view of the hypothesis
($f_1$),  $f(x,s)=as^-+bs^++g(x,s)$, $g(x,0)=0$, $\lim\limits_{|s|\to\infty}\frac{g(x,s)}s=0$ uniformly with respect to $x\in \mathbb{R}^N$, we obtain
\begin{equation*}
\frac {f^*_m(x,s)}s\to \begin{cases}
a+m & as  \quad s \to -\infty,\\
0 & as \quad s \to +\infty.
\end{cases}
\end{equation*}

From the hypothesis
($f_2$), $f(x,s)=a_0s^-+b_0s^++\tilde g(x,s)$, $\lim_{|s|\to 0_-}\frac{f(x,s)}s=a_0$,  and $\lim_{|s|\to 0_+}\frac{f(x,s)}s=b_0$ uniformly on $x\in \mathbb{R}^N$, then $g(x,u)=\tilde g(x,u)+(a_0-a)u^-+(b_0-b)u^+$,
we obtain
\begin{align}\label{117}
\lim_{s\to 0_-}\frac{f^*_m(x,s)}s&=\lim_{s\to 0_-}\frac{(a+m)s^-+g(x,s^-)}s \nonumber\\ \cr&=\lim_{s\to 0_-}\frac{(a+m)s^-+\tilde g(x,s^-)+(a_0-a)s^-+0\cdot s^+}s=a_0+m,
\end{align}

\begin{align}\label{118}
\lim_{s\to 0_+}\frac{f^*_m(x,s)}s&=\lim_{s\to 0_+}\frac{(a+m)s^-+g(x,s^-)}s\nonumber\\  &=\lim_{s\to 0_-}\frac{(a+m)s^-+\tilde g(x,s^-)+(a_0-a)s^-+0\cdot s^+}s=0,
\end{align}
uniformly on $x\in \mathbb{R}^N$.

Set $f^*_m(x,s)=(a_0+m)s^-+0\cdot s^++\tilde g(x,s^-)$, where  $\lim\limits_{s\to 0}\frac{\tilde g(x,s^-)}s=0$ uniformly with respect to $x\in \mathbb{R}^N$.

\subsection{The properties of spectrum of Schr\"{o}dinger operator  with  certain potentials }\label{3.2}

To show our consequence, we use the following well-known theorem from  spectral theory.

\begin{proposition}\label{prop31}(see Theorem XIII.48, \cite{18})

Let $G$ be an closed subset of $\mathbb{R}^N$  of measure zero and suppose that $V\in L^1_{loc}(\mathbb{R}^N\setminus G)$  with $V\geq 0$. Let $H=-\Delta +V$ be defined as the sum of quadratic forms. Then :

(a) If $\mathbb{R}^N\setminus G$ is connected and $\inf \sigma(H)$ is an eigenvalue, then it is a simple eigenvalue and the corresponding eigenfunction is strictly positive.\\

(b) If $\mathbb{R}^N\setminus G$ is disconnected, $V$ can be chosen so that $H$ has a degenerate ground state.

\end{proposition}

Suppose $V=V^++V^-$  is real-valued,  $V^+\in L^r_{loc}(\mathbb{R}^N)$ and $V^-$ is a K-R potential. Assume $V^-=\tilde V^*_1+\tilde V^*_2$, $ \tilde V^*_1\in L^p(\mathbb{R}^N)$, $\tilde V^*_2\in L^{\infty}(\mathbb{R}^N)$, $p=2$ if $N\leq 3$, $p> 2$ if $N=4$, $p>\frac N2$ if $N\geq 5$.\\

Set $V_n=V^++\tilde V^*_{1,n}+\tilde V^*_2$, where $\tilde V^*_{1,n}=\begin{cases}
\tilde V^*_1,& if  \quad \tilde V^*_1\geq -n,\\
-n,& if  \quad \tilde V^*_1< -n,
\end{cases}$.
It is obvious that $ V_n$ is bounded below, $V_n\geq V$ a.e. on $\mathbb{R}^N$ and
$(V_n- V)=(\tilde V^*_{1,n}-\tilde V^*_1)\in L^p(\mathbb {R}^N)$.

In the proofs that follow, we shall routinely take $N\geq 3$. The proofs for $N=1$ or $2$ are not more complicated. To begin, we show

\begin{lemma}
$(V_n- V)\to 0$ in $L^p(\mathbb {R}^N)$.
\end{lemma}
\textit{Proof.}
Notice that  $(\tilde V^*_{1,n}-\tilde V^*_1)\to 0$ in $L^p(\mathbb {R}^N)$ as $n\to \infty$, and this implies $(V_n- V)\to 0$ in $L^p(\mathbb {R}^N)$.

Indeed since  $\tilde V^*_{1,n}, \tilde V^*_1\in L^p(\mathbb {R}^N)$, $ \tilde V^*_{1,n}\to \tilde V^*_1 $   almost everywhere,  and  $\|\tilde V^*_{1,n}\|_{L^p}\to \|\tilde V^*_{1}\|_{L^p}$, by Br\'{e}zis-Lieb Lemma we have $(\tilde V^*_{1,n}-\tilde V^*_1)\to 0$ in $L^p(\mathbb {R}^N)$ as $n\to \infty$. \qed \vskip 5pt

Denote $A_m=-\Delta+V+m$, and $A^{(n)}_m=-\Delta +V_n+m$. They have the same domain $H^2_m(\mathbb{R}^N)=\{u\in L^2(\mathbb{R}^N): (-\Delta+V^+ )u\in  L^2(\mathbb{R}^N)\}$ and form domain $H^1_m(\mathbb{R}^N)=\{ u\in H^1(\mathbb{R}^N): \int_{\mathbb{R}^N}|\nabla u|^2+(V^++m) u^2< +\infty \}$(see Remark \ref{remk5} for  details).\\

Set
\begin{equation}
\lambda_1=\inf_{u\in H^1_m,\|u\|_{L^2}=1}\int_{\mathbb{R}^N}|\nabla u|^2+V u^2,
\end{equation}
\begin{equation}
\lambda_1^{(n)}=\inf_{u\in H^1_m,\|u\|_{L^2}=1}\int_{\mathbb{R}^N}|\nabla u|^2+V _n u^2.
\end{equation}

\begin{lemma}\label{lem36}
Under the above assumptions, (i) $\lambda_1^{(n)}\to \lambda_1$, (ii) $\lambda_1^{(n)}<\inf\sigma_{ess}(-\Delta+V_n)$, for $n$ large.
\end{lemma}
\textit{Proof.}
First we show $\lambda_1^{(n)}$ is bounded below.
Obviously,
for any $u\in H^1_m(\mathbb{R}^N)$,
\begin{equation}\label{1331}
\int_{\mathbb{R}^N}|\nabla u|^2+V _n u^2=\int_{\mathbb{R}^N}|\nabla u|^2+V u^2+\int_{\mathbb{R}^N}(V_n-V) u^2\geq \int_{\mathbb{R}^N}|\nabla u|^2+V u^2,
\end{equation}
we obtain
$\lambda_1^{(n)}\geq \lambda_1$.

Next we prove  $\lambda_1^{(n)}$ is bounded above.
As $\lambda_1\in \sigma_{dis}(A=-\Delta +V)$,
denote by $\phi_0\in H^2_m$ one of the corresponding eigenvector of $\lambda_1$, and $\|\phi_0\|_{L^2}=1$.
Notice $\|\phi_0\|^2_{_m}=(\lambda_1+m)\|\phi_0\|^2_{L^2}$.
Then
\begin{align}\label{1332}
\lambda_1^{(n)}&\leq \int_{\mathbb{R}^N}|\nabla \phi_0|^2+V_n \phi_0^2\nonumber\\
               &=   \int_{\mathbb{R}^N}|\nabla \phi_0|^2+(V_n-V)\phi_0^2+V\phi_0^2\nonumber\\
               &=  \int_{\mathbb{R}^N}|\nabla \phi_0|^2  +V\phi_0^2+     \int_{\mathbb{R}^N}(V_n-V)\phi_0^2\nonumber\\
               &\leq  \int_{\mathbb{R}^N}|\nabla \phi_0|^2  +V\phi_0^2+ \|V^-_n-V^-\|_{L^p}\|\phi_0\|^2_{L^{2q}}\nonumber\\
               &\leq \lambda_1+\varepsilon,
\end{align}
since $\|V^-_n-V^-\|_{L^p}\to 0$ and $\|\phi_0\|^2_{L^{2q}}\leq C_0\|\phi_0\|^2_{m}=C_0(\lambda_1+m)$,  $\frac1p+\frac1q=1$.

That's the precise statement (i) by combing (\ref{1331}) with (\ref{1332}).\\
Aiming at verifying (ii), we need the following proposition.

\begin{proposition}(see Theorem 14.11, \cite{10})
Let $V$ be a real-valued potential in the Kato-Rellich class, and let $H=-\Delta +V$ be the corresponding self-adjoint, semi-bounded  Schr$\ddot{o}$dinger operator with domain $H^2(\mathbb{R}^N)$. Then the bottom of the essential spectrum is given by
\begin{equation}\label{147}
\inf \sigma_{ess}(H)=\sup _{\kappa\subset\mathbb{R}^N}\inf _{\phi\neq 0,\phi \in C^{\infty}_0(\mathbb{R}^N\setminus \kappa)}\frac {\langle\phi ,H\phi\rangle}{\|\phi\|^2_{L^2}},
\end{equation}
where the supremum is over all compact subset $ \kappa \subset \mathbb{R}^N $.
\end{proposition}

The formula is true (see \cite{16}), if $H=-\Delta + a(x)$, where $a(x)$ is a real measurable function and satisfies
\begin{equation}\label{1111}
  a_0=\lim_{r\to \infty}(ess\inf_{x\in \mathbb{R}^N\setminus B_r(0)}a(x))>-\infty,
\end{equation}
and
\begin{equation}\label{222}
  \int_K|a(x)|^{\frac N2+\varepsilon}dx<\infty,
\end{equation}
for any bounded set $K\subset\mathbb{R}^N$ and some $\varepsilon>0$. In particular, if $V$ is real-valued,  $V=V^++V^-$, $ V^+\in L^r_{loc}(\mathbb{R}^N)$ and $V^-$ is a K-R potential, $V$ satisfies (\ref {1111}) and (\ref {222}), where $r=2$ if $N\leq 3$,  and $r>\frac N2$, if $N\geq 4$.
\begin{lemma}
Under the assumption above, $\lambda_1^{(n)}\notin \sigma_{ess}(-\Delta+V_n)$.
\end{lemma}

\textit{Proof.}
Claim that for $n$ large, $\lambda_1^{(n)}<\inf \sigma_{ess}(-\Delta+V_n)$.
Since $V_n\geq V$ a.e on $\mathbb{R}^N $, for $\forall \phi\in C^{\infty}_0(\mathbb{R}^N\setminus \kappa), \phi \neq 0$, we have

\begin{equation}\label{34}
\frac{\langle\phi ,(-\Delta+V_n)\phi\rangle_{L^2}}{\|\phi\|^2_{L^2}}\geq \frac{\langle\phi ,(-\Delta+V)\phi\rangle_{L^2}}{\|\phi\|^2_{L^2}}.
\end{equation}
Since
\begin{equation}\label{es}
\inf \sigma_{ess}(-\Delta+V_n)=\sup _{\kappa\subset\mathbb{R}^N}\inf _{\phi\neq 0,\phi \in C^{\infty}_0(\mathbb{R}^N\setminus \kappa)}\frac{\langle\phi ,(-\Delta+V_n)\phi\rangle_{L^2}}{\|\phi\|^2_{L^2}},
\end{equation}

\begin{equation}\label{se}
\inf \sigma_{ess}(-\Delta+V)=\sigma_0=\sup _{\kappa\subset\mathbb{R}^N}\inf _{\phi\neq 0,\phi \in C^{\infty}_0(\mathbb{R}^N\setminus \kappa)}\frac{\langle\phi ,(-\Delta+V)\phi\rangle_{L^2}}{\|\phi\|^2_{L^2}}.
\end{equation}
combining (\ref{34}), (\ref{es}) and (\ref{se}) we get
\begin{equation}\label{ess}
\inf \sigma_{ess}(-\Delta+V_n)\geq \inf \sigma_{ess}(-\Delta+V)=\sigma_0>\lambda_1.
\end{equation}
We now go on the proof procedure of Lemma \ref{lem36}.
It deduces from (\ref{ess}) that for $n$ large, $\lambda_1^{(n)}<\inf \sigma_{ess}(-\Delta+V_n)$ via $\lambda_1^{(n)}\to \lambda_1$.
Therefore we may apply Proposition \ref{prop31} to obtain that $\lambda_1^{(n)}$ is a simple eigenvalue of $-\Delta+V_n$ for all $n$  sufficiently large and the corresponding eigenfunction is strictly positive. The proof is complete.\qed \vskip 5pt
\begin{corollary}
Similarly,   since $\lambda_1$ and $\lambda_1^{(n)}$  are simple, we obtain, for $n$ large there holds $\lambda_2^{(n)}\geq \lambda_2$.
\end{corollary}

\begin{lemma}\label{LM1}
Assume that $(V_1)$ and  $(V_2)$ hold. Then  $\lambda_1$ is an eigenvalue of $-\Delta+V $ associated to an eigenvector    $\varphi_1$ which is non-negative.  Moreover, for any $a>\inf\sigma(-\Delta+V )$ there holds $(a,0)\bigcup (0,a) \notin \Sigma(-\Delta+V )$.
\end{lemma}

\textit{Proof.}  Let $\lambda_1^{(n)}$ denote the simple eigenvalue of $-\Delta+V_n$ and the corresponding eigenfunction $\varphi_1^{(n)}$ is strictly positive.
Set $\|\varphi_1^{(n)}\|^2_{L^2}=1$. Obviously, $\{\varphi_1^{(n)}\}$ is bounded in $H^1_m(R^N)$. Up to a subsequence, $\varphi_1^{(n)}\rightharpoonup w$ in $H^1_m(\mathbb{R}^N)$.

Define $I(u)=\frac 12 \int_{\mathbb{R}^N}|\nabla u|^2  +V(x) u^2 -\frac 12\lambda_1\int_{\mathbb{R}^N}u^2 dx$, we have $I^{\prime}(\varphi_1^{(n)})\to 0$, as $n \to \infty$. In fact,

\begin{align*}
\langle I^{\prime}(\varphi_1^{(n)}),v\rangle_{m}&= \int_{\mathbb{R}^N}\nabla \varphi_1^{(n)}\nabla v+V_n \varphi_1^{(n)} v-\lambda_1 \varphi_1^{(n)}v+(V-V_n)\varphi_1^{(n)} v\nonumber\\
                                    &=\int_{\mathbb{R}^N}(\lambda^{(n)}_1-\lambda_1)\varphi^{(n)}_1v+(V-V_n)v\varphi_1^{(n)}\nonumber\\
                                    &\leq (\lambda^{(n)}_1-\lambda_1)\|\varphi_1^{(n)}\|_{L^{2}} \|v\|_{L^{2}}+ \|V^-_n-V^-\|_{L^p}\|v\|_{L^{2q}}\| \varphi_1^{(n)}\|_{L^{2q}}\nonumber\\
                                    &\to 0,
\end{align*}
since $\lambda^{(n)}_1\to \lambda_1, \|V^-_n-V^-\|_{L^p}\to 0$ and $\|\varphi_1^{(n)}\|_{L^{2q}}\leq C_0\|\varphi_1^{(n)}\|_{m}\leq C_0(\lambda_1+m+1)$,  where $\frac1p+\frac1q=1$.

$\bf{Step 1.}$
We show that $w$ is an eigenvector of  $-\Delta+V$  associated to the eigenvalue  $\lambda_1$.\\
It is sufficient to show that for any $v\in C_0^{\infty}(\mathbb{R}^N)$
\begin{align}\label{3356}
\int_{\mathbb{R}^N}\nabla w\nabla v  +V(x) wv dx=\lambda_1\int_{\mathbb{R}^N}wv dx.
\end{align}
Let $v\in C_0^{\infty}(\mathbb{R}^N)$ be arbitrary but fixed. Since ${\varphi_1^{(n)}}\rightharpoonup w$  in $H^1_m(\mathbb{R}^N)$, we have
\begin{align}\label{3357}
\int_{\mathbb{R}^N}\nabla \varphi_1^{(n)}
\nabla v  +(V(x)+m) \varphi_1^{(n)}v dx\to \int_{\mathbb{R}^N}\nabla w
\nabla v  +(V(x)+m) wv dx.
\end{align}
By the choice of $ \varphi_1^{(n)}$, we obtain
\begin{align}\label{3358}
\int_{\mathbb{R}^N}\nabla \varphi_1^{(n)}
\nabla v  +(V(x)+m) \varphi_1^{(n)}v dx &=   \int_{\mathbb{R}^N} \nabla \varphi_1^{(n)}\nabla v +(V_n +m) \varphi_1^{(n)}v+(V-V_n)u\varphi_1^{(n)}\nonumber\\
               &= \lambda_1^{(n)} \int_{\mathbb{R}^N}v  \varphi_1^{(n)}+ \int_{\mathbb{R}^N}(V-V_n)v\varphi_1^{(n)}\nonumber\\
               &\leq\lambda_1^{(n)} \int_{\mathbb{R}^N}v  \varphi_1^{(n)}+ \|V^-_n-V^-\|_{L^p}\|v\|_{L^{2q}}\| \varphi_1^{(n)}\|_{L^{2q}}\nonumber\\
               &= (\lambda_1^{(n)}+m) \int_{\mathbb{R}^N}v  \varphi_1^{(n)}+\varepsilon_n \nonumber\\
               &\to(\lambda_1+m) \int_{supp v}v  w.
\end{align}
Combining (\ref{3357}) with (\ref{3358}) we get (\ref{3356}). Next we prove $w\neq 0$.
Seeking a contradiction, suppose that $w=0$.
Let $P_{N_1}$ denote the orthogonal projection of $-\Delta+V $ onto $N_1=span\{\varphi_1\}$ in $L^2(\mathbb{R}^N)$, where $\varphi_1$ is the eigenvector of  $-\Delta+V$ associated to the eigenvalue  $\lambda_1$,
and $P_{N_1^{\perp}}$ denotes the orthogonal projection onto $N_1^{\perp}$ in $L^2(\mathbb{R}^N)$.

It is easy to check that
\begin{align*}
&\int_{\mathbb{R}^N} |\nabla P_{N_1}\varphi_1^{(n)}|^2+(V(x)+m)(P_{N_1}\varphi_1^{(n)})^2\nonumber\\
&+\int_{\mathbb{R}^N} |\nabla P_{N_1^{\perp}}\varphi_1^{(n)}|^2+(V(x)+m)(P_{N_1^{\perp}}\varphi_1^{(n)})^2 \nonumber\\ &\leq  (\lambda_1+m) \int_{\mathbb{R}^N}(\varphi_1^{(n)})^2dx+\varepsilon_n.
\end{align*}

Since $P_{N_1}$ has $1$-dimensional and $w=0$, we get $P_{N_1}\varphi_1^{(n)}\to 0$ in $H^1_m(\mathbb{R}^N)$.  Thus
\begin{align*}
&\int_{\mathbb{R}^N} |\nabla P_{N_1^{\perp}}\varphi_1^{(n)}|^2
+(V(x)+m)(P_{N_1^{\perp}}\varphi_1^{(n)})^2\nonumber\\
&\leq  (\lambda_1+m+\varepsilon_n)  \int_{\mathbb{R}^N}(P_{N_1^{\perp}}\varphi_1^{(n)})^2dx\nonumber\\
&\leq\frac{\lambda_1+m +\varepsilon_n}{\lambda_2+m}\int_{\mathbb{R}^N} |\nabla P_{N_1^{\perp}}\varphi_1^{(n)}|^2
+(V(x)+m)(P_{N_1^{\perp}}\varphi_1^{(n)})^2,
\end{align*}
namely,
\begin{equation*}
\| P_{N_1^{\perp}}\varphi_1^{(n)}\|_m^2\leq \frac{\lambda_1+m+\varepsilon_n}{\lambda_2+m} \| P_{N_1^{\perp}}\varphi_1^{(n)}\|_m^2,
\end{equation*}
which is impossible.

Since
\[\langle I^{\prime}(\varphi_1^{(n)})-I^{\prime}(w),\varphi_1^{(n)}-w\rangle_{m}\to 0,\]
we have
\begin{align*}
&\int_{\mathbb{R}^N} |\nabla P_{N_1}(\varphi_1^{(n)}-w)|^2+(V(x)+m)(P_{N_1}(\varphi_1^{(n)}-w))^2\nonumber\\
&+\int_{\mathbb{R}^N} |\nabla P_{N_1^{\perp}}(\varphi_1^{(n)}-w)|^2+(V(x)+m)(P_{N_1^{\perp}}(\varphi_1^{(n)}-w))^2\nonumber\\
&\leq  (\lambda_1+m) \int_{\mathbb{R}^N}(\varphi_1^{(n)}-w)^2dx+\varepsilon_n.
\end{align*}
Since $P_{N_1}$ has 1-dimensional and $w\neq 0$,  we get
\begin{equation*}
\| P_{N_1^{\perp}}(\varphi_1^{(n)}-w)\|_m^2\leq \frac{\lambda_1+m+\varepsilon_n}{\lambda_2+m} \| P_{N_1^{\perp}}(\varphi_1^{(n)}-w)\|_m^2,
\end{equation*}
which means $\varphi_1^{(n)}\to w$ in $H^1_m(\mathbb{R}^N)$.

Since $\varphi_1^{(n)}$  is strictly positive, we have that
the weak limit $w$ is a non-negative eigenvector of $-\Delta+V$ associated to the eigenvalue $\lambda_1$. This ends Step 1.

$\bf{Step  2.}$ When $a>\inf\sigma(-\Delta+V )$, the operator $-\Delta+V $  has no non-negative eigenvector associated to the eigenvalue $a$.
Seeking a contradiction, suppose that $u\in H^1_m(\mathbb{R}^N)$ is non-negative and satisfies
\[ -\Delta u+V u=au \quad \text{in}\quad R^N.
\]
First, we fix a constant $l$ such that
\[ \lambda_1=\inf\sigma(-\Delta +V)<l<a.
\]
By the variational characterization of  $\inf\sigma(-\Delta+V_n)$,  there exists $N_0$ satisfying   $\lambda_1^{(N_0)}<l$, since $\lambda_1^{(n)}\to \lambda_1$.
Suppose that $\varphi_1^{(N_0)}$   is the eigenvector of  $-\Delta+V_{N_0}$ associated to the eigenvalue $\lambda_1^{(N_0)}$. Since $V-V_{N_0}\leq 0$ and  $\varphi_1^{(N_0)}>0$ then

\begin{align*}
a\langle u,\varphi_1^{(N_0)}\rangle_{L^2} &= \langle (-\Delta+V)u ,\varphi_1^{(N_0)}\rangle_{L^2}\nonumber\\
               &=   \int_{\mathbb{R}^N}\nabla u\nabla \varphi_1^{(N_0)} +V u\varphi_1^{(N_0)}\nonumber\\
               &=   \int_{\mathbb{R}^N}\nabla u\nabla \varphi_1^{(N_0)} +V_{N_0} u\varphi_1^{(N_0)}+(V-V_{N_0})u\varphi_1^{(N_0)}\nonumber\\
               &= \langle u, (-\Delta+V_{N_0})\varphi_1^{(N_0)}\rangle_{L^2} + \int_{\mathbb{R}^N}(V-V_{N_0})u\varphi_1^{(N_0)}\nonumber\\
               &\leq \langle u, (-\Delta+V_{N_0})\varphi_1^{(N_0)} \rangle_{L^2}\nonumber\\
               &=\lambda_1^{(N_0)} \int_{\mathbb{R}^N}u  \varphi_1^{(N_0)}\nonumber\\
               &\leq l\langle u,\varphi_1^{(N_0)}\rangle_{L^2},
\end{align*}
where $\langle \cdot,\cdot\rangle_{L^2}$ denotes the scalar product of $L^2(R^N)$. But since $u\geq 0$ and $\varphi_1^{(N_0)}>0$, we have $\langle u,\varphi_1^{(N_0)}\rangle_{L^2} >0$, and thus the above calculation shows that $a\leq l $ in contradiction with the choice of $l$. This completes Step 2.

Combining Step 1 and Step 2, Lemma \ref{LM1} is proved.\qed \vskip 5pt
An immediate consequence of Lemma \ref{LM1} is:
\begin{corollary}\label{cor8}
If $a>\inf\sigma(-\Delta+V)$, then
$(a+m,0)\notin \sum(-\Delta+V+m)$.
\end{corollary}

Next, we will use the computations of critical groups of $J$ at infinity  and zero with  topological analysis of level sets of $J$ to infer the existence of nontrivial solutions of equation (\ref{fix}), via Morse theory.
An essential step in  these arguments is the $(PS)$  condition for the functional.  We will use the following version of $(PS)$ condition, derived in \cite{12}.

\begin{theorem}\label{them2}
Let $V$ be a real potential as above, $\sigma_{dis}(-\Delta+V)\ne \varnothing$, and $\inf\sigma(-\Delta+V)\in \sigma_{dis}(-\Delta+V)$. Under the hypothesis $(f_1)$, if $a \neq \lambda_1$, then $J$ satisfies the $(PS)$ condition on $H^1_m(\mathbb{R}^N)$.
\end{theorem}

\textit{Proof.}
This is actually a slight variation of Theorem 3.4  in \cite{12}  for which the same proof applies, so we omit it. \qed \vskip 5pt

\subsection{Computation of $C_*(J,\infty)$}\label{3.4}

At this point, we will reduce  the computation of  critical groups   of  $J$ at  infinity to   compute critical groups of corresponding  homogeneous functional at zero, inspired by the paper \cite{12}.

We need some notions and known results on Morse theory.
Let H be a Hilbert space and $f\in C^1(H,\mathbb {R})$. $K_f=\{u\in H:f^{\prime}(u)=0\}$ is the set of critical points of $f$ on $H$, and assume that $u_0$ is an isolated critical point of $f$ on $H$, $c=f(u_0)$. The $q$th critical group with coefficient group $G$ of $f$ at $u_0$ is defined by
\begin{equation*}
C_q(f,u_0)=H_q(f^c\cap U,(f^c\setminus \{u_0\})\cap U,G),
\end{equation*}
where $f^c=\{u\in H:f(u)\leq c\}$ denotes the sub-level sets, $U$ is any neighborhood of $u_0$ such that $u_0$ is the only critical point of $f$ in $f^c\cap U$.

Let $H$ be a Hilbert space and $f\in C^1(H, \mathbb{R})$. Recall the concept of critical groups of $f$ at infinity introduced by \cite{1}.

\begin{definition}
Suppose $f(K_f)$ is strictly bounded from below by $c\in  \mathbb{R}$ and that $f$ satisfies the $(PS)$ condition. Then $ C_q(f,\infty)=H_q(H,f^c)$, $q\geq 0$, is the $q$th critical group of $f$ at infinity. It is independent of the choice of $c$ with the $(PS)$ condition.
\end{definition}
To obtain one nontrivial solution, we must determine the critical groups of $0$ and $\infty$.
Our main result concerning the computation of critical groups at infinity reads:
\begin{theorem}\label{themJ}
Suppose $(V_1), (V_2)$ and  $(f_1)$ hold. If $\sigma_0>a>\lambda_1$, then
\begin{equation}\label{1234}
   C_q(J,\infty)=0, q=0,1,2\cdots.
\end{equation}

\end{theorem}
The proof of Theorem \ref{themJ} requires some preliminary lemmas.

Define
\begin{equation}\label{Jnt}
  J_{t}(u)=\frac12\int_{\mathbb{R}^N} |\nabla u|^2 +(V+m)u^2 dx-\frac {a+m}2\int_{\mathbb{R}^N} (u^-)^2 -t \int_{\mathbb{R}^N} G^*(x,u).
\end{equation}
$t\in[0,1], G^*(x,u)=\int_0^{u^-}g(x,s)ds$, and consider the following equation:
\begin{equation*}
\begin{cases}
-\Delta u+(V+m)u=(a+m)u^-+tg(x,u^-),& x\in \mathbb {R}^N,\\
u(x)\to 0,& as |x|\to +\infty,
\end{cases}
\end{equation*}
and define
\begin{equation}\label{22}
I_{(a+m,0)}(u) =\frac12\int_{\mathbb{R}^N} |\nabla u|^2+(V(x)+m)u^2 dx-\frac {1}2\int_{\mathbb{R}^N} (a+m)(u^-)^2.
\end{equation}

\begin{lemma}
Under the above hypotheses, there exists $M>0$, s.t.
\begin{equation}\label{63}
\sup_{t\in[0,1]}\sup_{u\in K_{J_{t}}}\|u\|_m<M.
\end{equation}
\end{lemma}
\textit{Proof.} This is actually a slight variation of Theorem 3.4 in \cite{12}  for which the same proof applies, so we omit it.
\qed \vskip 5pt

\begin{lemma}
Under the   above hypotheses,
\begin{equation}\label{64}
\delta =\inf_{t\in[0,1],u\in H^1_m(\mathbb{R}^N )\setminus B_m(0,M)}\frac {\|J^{\prime}_{t}(u)\|_m}{\|u\|_m}>0,
\end{equation}
where
$J^{\prime}_{t}(u)=u-A_m^{-1}[(a+m)u^- +tg(x,u^-)]$, $B_m(u,M)=\{\omega \in H^1_m (\mathbb{R}^N ):\|\omega-u\|_m<M\}$, where $M$ is given by (\ref{63}).

\end{lemma}
\textit{Proof.} By way of negation, there exists $t_k\in [0,1]$, $\{u_{k}\}_{k=1}^{\infty}\in H^1_m(\mathbb{R}^N )\setminus B_m(0,M)$, s.t.
\begin{equation}\label{65}
\frac {\|J^{\prime}_{t_k}(u_{k})\|_m}{\|u_{k}\|_m}\to 0, k\to +\infty.
\end{equation}
Repeating an argument analogous to step two of the proof of Theorem 3.4 in \cite{12}, it follows that $\{u_{k}\}_{k=1}^{+\infty}$ is bounded in $ H^1_m(\mathbb{R}^N)$, and $\|J^{\prime}_{t_k}(u)\|_m\to 0$. Assume $t_k\to t_0$, $u_{k}\rightharpoonup u_{0}$ in $ H^1_m(\mathbb{R}^N)$, as $k\to \infty$, we get
\begin{equation*}
\int_{\mathbb{R}^N} \nabla u_{0}\nabla \varphi+ (V(x)+m)u_{0}\varphi=\int_{\mathbb{R}^N} (a+m)u^-_{0}\varphi +t_0g(x,u^-_{0})\varphi,\nonumber
\end{equation*}
for any $ \varphi\in C^{\infty}_0(\mathbb{R}^N)$.\\
Based on the density of $ C^{\infty}_0(\mathbb{R}^N)$ in $H^1_m(\mathbb{R}^N)$, we yield
\begin{equation*}
\int_{\mathbb{R}^N} \nabla u_0\nabla \varphi+ (V(x)+m)u_{0}\varphi=\int_{\mathbb{R}^N} (a+m)u^-_{0}\varphi+t_0g(x,u^-_{0})\varphi , \forall \varphi\in H^1_m(\mathbb{R}^N).
\end{equation*}
and this indicates that $u_{0}$ is a weak solution of (\ref{Jnt}) for $t=t_0$, i.e.,
\begin{equation*}
J^{\prime}_{t_0}(u_{0})=u_{0}-A_m^{-1}[(a+m)u_{n,0}^-+t_0g(x,u^-_{0})]=0.
\end{equation*}
Notice that
\begin{equation*}
J^{\prime}_{t_k}(u_{k})=u_{k}-A_m^{-1}[(a+m)u_{k}^-+t_k g(x,u^-_{k})]\to  0.
\end{equation*}
Fully imitating the trick employed by the proof of Theorem 3.4 in \cite{12}, we get $u_{k}\to u_{0}$ in $H^1_m(\mathbb{R}^N)$, so $u_{0}\in K_{J_{t_0}}\setminus B_m(0,M)$, which contradicts (\ref{63}). The proof is complete.\qed \vskip 5pt

\begin{lemma}
Under the hypothesis $(f_1)$,  there exists $\tilde C<0$,  for any $ t\in[0,1]$,  $ u\in H^1_m(\mathbb{R}^N)$, if $J_{t}(u) \leq \tilde C$, then $\|u\|_m \geq M$.
\end{lemma}
\textit{Proof.}  Assume the contrary.  Take $C_k\to -\infty$ as $k\to +\infty$, then $\exists t_k \in[0,1]$ and $u_k\in E$, $J_{t_k}(u_k) \leq C_k$, $\|u_k\|_m <M$. We just treat the case $a+\beta >\lambda_1$. Observe that

\begin{align}\label{610}
J_{t_k}(u_k)&\geq\frac12 \| u_k\|_m^2-\frac {(a+m+\beta)}2\int_{\mathbb{R}^N} u_k^2\nonumber\\
 &\geq\frac12 \| u_k\|_m^2-\frac{(a+m+\beta)}{2(\lambda_1+m)}\| u_k\|_m^2\nonumber\\
 &\geq -\frac{(a+\beta-\lambda_1)}{2(\lambda_1+m)}\| u_k\|_m^2\nonumber\\
 &\geq -\frac{(a+\beta-\lambda_1)}{2(\lambda_1+m)}M^2\nonumber\\
\end{align}
as $k\to +\infty$. A contradiction!\qed \vskip 5pt

\begin{theorem}\label{them16}
Assume  $(V_1)$, $(V_2)$, $(f_1)$ hold. If $(a+m,0)\notin \Sigma (-\Delta +V+m)$, and $ a<\sigma_0,$ then
\begin{equation}\label{635}
 C_q(J,\infty)\cong C_q(J_{t},\infty)=C_q(I_{(a+m,0)},0), \forall t\in[0,1].
\end{equation}
\end{theorem}
\textit{Proof.}
This is actually a slight variation of Theorem 6.2 in \cite{12}  for which the same proof applies, so we omit it.

\subsubsection{Computation of $C_q(I_{(a+m,0)},0)$}

\begin{lemma}\label{lem17} Assume $(f_1)$ holds. If $a>\lambda_1$, then
\begin{equation}\label{i}
 C_q(I_{(a+m,0)},0)=0, q=0,1,2\cdots.
\end{equation}
\end{lemma}
\textit{Proof.}
Since $\sigma_0>a>\lambda_1$, we have $(a+m,0)\notin \Sigma(-\Delta+V+m)$, then it is easy to check that $C_q(I_{(a+m,0)},0)=C_q(I_{(\hat a+m,0)},0)$ for any $\hat a\in (\lambda_1,\sigma_0)$.

It suffices to show that if $\lambda_2>a>\lambda_1$, then
\begin{equation}\label{i}
 C_q(I_{(a+m,0)},0)=0,q=0,1,2\cdots.
\end{equation}

Next we prove
\begin{lemma}\label{lem8} Assume $(f_1)$ holds. If $\lambda_2>a>\lambda_1$, then
\begin{equation}\label{i}
 C_q(I_{(a+m,0)},0)=0,q=0,1,2\cdots.
\end{equation}
\end{lemma}
\textit{Proof.} Since $\lambda_1=\inf\sigma(-\Delta+V)$ is simple, we denote by $ \varphi_1\geq 0$  the corresponding eigenvector of $\lambda_1$, and set
$\|\varphi_1\|_{L^2}=1$.

Denote by $P_0$ the orthogonal projection of $H^1_m(\mathbb {R}^N)$ onto the space $E_1=\text{span} \{\varphi_1 \}$, and $P_{\perp}$ the projection onto the space  $E_1^{\perp}$ in $L^2(\mathbb {R}^N)$.
For any $u\in H^1_m(\mathbb {R}^N)$, write $u=u_1+u_*$, where  $u_1=\alpha\varphi_1$,
$P_0u=\alpha\varphi_1$, $P_{\perp}u=u_*$.

Denote
\begin{equation}\label{it}
I_t(u)=(1-t)I_{(a+m,0)}(u)+t(\frac{ \text{sign } \alpha\|\alpha\varphi_1\|_m^2+\|u_{*}\|_m^2}2).
\end{equation}
Note that $\int_{\mathbb{R}^N} \varphi u_*=0$ and $\langle \varphi,u_*\rangle_m=0$.

Claim for any $u\neq 0$, $\forall t\in [0,1]$, $ \exists v\neq 0$ such that $\langle I^{\prime}_t(u),v\rangle \neq 0$, i.e $I^{\prime}_t(u)\neq 0$.
We consider two cases:

(1) $u^{+}\ne 0$.

If $\alpha\leq0$, then
take $v=u^+$.

\begin{eqnarray*}
 & &\langle I^{\prime}_t(u) , u^+\rangle_m\nonumber\\ &=&(1-t)[\langle u,u^+\rangle_m-(a+m)\int_{\mathbb{R}^N} u^-u^+]+t[\langle u,u^+\rangle_m-\langle2\alpha\varphi_1,u^+\rangle_m]\nonumber\\
                      &=&\langle u,u^+\rangle_m-(1-t)(a+m)\int_{\mathbb{R}^N} u^-u^+-t\langle 2\alpha\varphi,u^+\rangle_m\nonumber\\
                      &\geq& \|u^+\|_m^2>0.
\end{eqnarray*}

If $\alpha>0$. Set $v=u^+$. Then

\begin{align*}
 \langle I^{\prime}_t(u) , u^+\rangle_m &=(1-t)[\langle u,u^+\rangle_m-(a+m)\int_{\mathbb{R}^N} u^-u^+]+t\langle u,u^+\rangle_m\nonumber\\
                      &=\langle u,u^+\rangle_m=\|u^+\|_m^2\nonumber\\
                      &>0.
\end{align*}

(2) $u^+=0$, i.e $u=u^-\leq 0$, then we have $\alpha=\langle u,\varphi_1\rangle_m\leq 0$.

(i) $u_*\neq0$. Take $v=u_*$.
\begin{eqnarray*}
 & &\langle I^{\prime}_t(u) , u_*\rangle_m\nonumber\\ &=&(1-t)[\langle u,u_*\rangle_m-(a+m)\int_{\mathbb{R}^N} u^-u_*]+t[\langle u,u_*\rangle_m-\langle 2\alpha\varphi_1,u_*\rangle_m]\nonumber\\
                      &=&\langle u,u_*\rangle_m-(1-t)(a+m)\int_{\mathbb{R}^N} uu_*\nonumber\\
                      &=& \|u_*\|_m^2-(1-t)(a+m)\int_{\mathbb{R}^N} (\alpha\varphi_1+u_*)u_*\nonumber\\
                      &=&\|u_*\|_m^2-(1-t)(a+m)\int_{\mathbb{R}^N} u^2_*\nonumber\\
                      &\geq&(1-\frac{a+m}{\lambda_2+m})\|u_*\|_m^2\nonumber\\
                      &>&0
\end{eqnarray*}

(ii) $u_*=0$. Then $u=\alpha \varphi_1$, $\alpha< 0$. Since $\lambda_1< a$, take $v=\alpha \varphi_1$, we have

\begin{eqnarray*}
& &\langle I^{\prime}(u),\varphi_1\rangle_m\nonumber\\
 &=&(1-t)\left( \langle\alpha \varphi_1,\varphi_1\rangle_m-(a+m)\int_{\mathbb{R}^N} u\varphi_1\right) +t \text{sign } \alpha\langle \alpha\varphi_1,\varphi_1\rangle_m \nonumber\\
                                  &=&(1-t)\left(\int_{\mathbb{R}^N}\alpha(\lambda_1+m)  \varphi_1^2-(a+m)\int_{\mathbb{R}^N}\alpha \varphi_1^2\right)+t\alpha \text{sign } \alpha\|\varphi_1\|_m^2  \nonumber\\
                                  &=&(1-t)\alpha (\lambda_1-a)\int_{\mathbb{R}^N}  \varphi_1^2+t\alpha \text{sign } \alpha\|\varphi_1\|_m^2  \nonumber\\
                                  &>&0.
\end{eqnarray*}

We proved that $0$ is the only critical point of $I_t$ in $B_r(0)$ if $r$ is sufficiently small. So
\begin{equation}\label{i01}
  C_q(I_t,0)=C_q(I_0,0)=C_q(I_1,0).
\end{equation}
Notice that $C_q(I_0,0)=C_q(I_{(a+m,0)},0)$,
we remain to compute the critical groups of $I_1$ at $0$.

Since
$I_1=\frac {(\text{sign } \alpha \|\alpha \varphi_1\|_m^2+\|u_*\|_m^2)}2$, on the unit ball $B_1$  of $H^1_m(\mathbb {R}^N)$ centered at $0$, it is easily seen that

\begin{equation}\label{0}
B_1\cap I_1^0=\{u\in H^1_m(\mathbb {R}^N) :\|u\|_m\leq 1,\alpha\leq 0 ,\text{and} \|\alpha \varphi_1\|_m\geq\|u_*\|_m\}.
\end{equation}\\
Define
\begin{equation}\label{t}
\eta(t,u)=\alpha\varphi_1+tu_*,\quad \forall(t,u)\in [0,1]\times(B_1\cap I_1^0).
\end{equation}

It is a strong deformation retract from $(B_1\cap I_1^0,B_1\cap (I_1^0\setminus\{0\}))$ to $(B_1\cap E_1^-,B_1\cap (E_1^-\setminus \{0\}))$, where $E_1= \text{span} \{\varphi_1 \}$ and $E_1^-= \{u\in E_1: u\leq 0 \}$. Thus
\begin{equation*}\begin{split}
C_q(I_1,0)&\cong H_q(I^0_1\cap B_1,(I^0_1\setminus \{0\})\cap B_1), \\
           &\cong H_q(B_1\cap E_1^-,B_1\cap (E_1^-\setminus \{0\}))\\
            &\cong H_q(X,Y).
\end{split}\end{equation*}
where $X=[-1,0],Y=[-1,0)$.
Notice that $X$ and $Y$ are contractible.
It is easy to check that $H_q(X,Y)=0$, $\forall q=0,1,2,\cdots$.
According to (\ref{i01}), we obtain $C_q(I_{(a+m,0)},0)=0$, $\forall q=0,1,2,\cdots $. The proof is complete.\qed \vskip 5pt

\textit{Proof of Theorem \ref{themJ}. }  By hypothesis $\sigma_0>a>\lambda_1$, and hence via Corollary \ref{cor8} we have $(a+m,0)\notin \Sigma(-\Delta+V+m)$.
Combining  Theorem \ref{them16} and Lemma \ref{lem17},  Theorem \ref{themJ} is proved.\qed \vskip 5pt

\subsection{Computation of $C_*(J,0)$}\label{3.5}
Set $f^*_m(x,s)=(a_0+m)s^-+\tilde g(x,s^-)$, $\tilde g(x,s^-)=g(x,s^-)-(a_0-a)s^--0\cdot s^+$, $g(x,0)=0$ and assume $\lim\limits_{s\to 0}\frac {\tilde g(x,s)}s=0$, $\lim\limits_{s\to \infty}\frac { g(x,s)}s=0$, uniformly with respect to $x \in \mathbb{R^N}$. To show our consequence, we need the following stronger hypothesis:

$(f^*_2) \max _{t\in [0,1]}a_t\leq \max _{t\in [0,1]}(a_t+t \beta) <\sigma_0,a_t=(1-t)a_0+ta$.\\

\begin{remark}
Note that the combination of $(f_1)$ and $a_0<\sigma_0$, alludes to $(f^*_2)$.
\end{remark}

Define
\begin{equation}\label{Jnt0}
\tilde J_{t}(u)=\frac12\int_{\mathbb{R}^N} |\nabla u|^2 +(V+m)u^2 dx-\frac {a_0+m}2\int_{\mathbb{R}^N} (u^-)^2 -t \int_{\mathbb{R}^N} \tilde G(x,u).
\end{equation}
$t\in[0,1], \tilde G(x,u)=\int_0^{u^-}\tilde g(x,s)ds$, and consider the following equation:
\begin{equation*}
\begin{cases}
-\Delta u+(V+m)u=(a_0+m)u^-+t\tilde g(x,u^-),& x\in \mathbb {R}^N,\\
u(x)\to 0,& as |x|\to +\infty,
\end{cases}
\end{equation*}
and define
\begin{equation}\label{220}
I_{(a_0+m,0)}(u) =\frac12\int_{\mathbb{R}^N} |\nabla u|^2+(V(x)+m)u^2 dx-\frac {1}2\int_{\mathbb{R}^N} (a_0+m)(u^-)^2.
\end{equation}
Our main result concerning the computation of critical groups at infinity reads:

\begin{theorem}\label{them5}
Suppose $(V_1), (V_2)$ and  $(f_2)$ hold. If $a_0>\lambda_1$, then
\begin{equation}\label{4}
   C_q(J,0)\cong C_q(\tilde J_t,0)\cong C_q(I_{(a_0+m,0)},0),\forall t\in [0,1].
\end{equation}

\end{theorem}

\textit{Proof.}
The proof is identical to that of Lemma 6.7 of \cite{12} and will be omitted.

\begin{lemma}\label{lem21}
$C_q(J,0)\cong \delta_{q0} \mathbb{Z}$.
\end{lemma}
\textit{Proof.}
It suffices to prove that $C_q(I_{(a_0+m,0)},0)\cong \delta_{q0} \mathbb{Z}$.  Assume
$\lambda_1 > a_0$, then
for any $ u\in H^1_m(\mathbb {R}^N)\setminus \{0\}$, we have

\begin{align}
I_{(a_0+m,0)}(u)&=\frac 12\int_{\mathbb{R}^N} |\nabla u|^2+(V(x)+m)u^2-(a_0+m)(u^-)^2dx\nonumber\\
                       &\geq \int_{\mathbb{R}^N}( \lambda_1+m) u^2-(a_0+m)(u^-)^2dx\nonumber\\
                       &\geq ( \lambda_1-a_0)\int_{\mathbb{R}^N} u^2\nonumber \\
                       &>0,
\end{align}
i.e., $ 0$ is a minimum of $I_{(a_0+m,0)}(u)$, then $C_q(I_{(a_0+m,0)},0)\cong \delta_{q0} \mathbb{Z}$.
The assertion comes from Theorem \ref{them5}  immediately. \qed\vskip 5pt

\section{Proof of Theorem \ref{them1}}\label{sec4}

Denote by $H_k(X,Y)$ the $k$th relative homology group with integer
coefficients. For $v$ an isolated critical point with $J(v)=c$ and $U$ a neighborhood of $v$, set
\[ C_k(J,v)=H_k(J^c\cap U,(J^c\cap U)\setminus \{v\}), k=0,1,\cdots,
\]
the $k$th critical group at the critical point $v$. If $J$ has only finitely many critical points $u_1,\cdots,u_n$ with $c_1<J(u_j)<c_2,$ we define the Morse numbers of the pair $(J^{c_2},J^{c_1})$ by
\[M_k=M_k(J^{c_2},J^{c_1})=\Sigma_{j=1}^{n} \dim C_k(J,u_j).
\]
Then, if $\beta_k=\dim H_k(J^{c_2},J^{c_1})$ are the Betti numbers of the pair $(J^{c_2},J^{c_1})$, the Morse inequalities require that:
\begin{equation}\label{morse}
\Sigma_{k=0}^{\infty}(-1)^k M_k=\Sigma_{k=0}^{\infty}(-1)^k \beta_k
\end{equation}
For a derivation of (\ref{morse}) and other facts from infinite-dimensional Morse theory see Chang \cite{6}.

\textit{Proof of Theorem \ref{them1}. }
To prove Theorem \ref{them1}, we argue by contradiction. Assume that $0$ is the only critical point of $J$.
According to lemma \ref{lem21}, $C_q(J,0)= \delta_{q0} \mathbb{Z}$,
hence the Morse numbers  are
\[M_q=\dim C_q(J,0)=\delta_{q0}, q=0,1,2\cdots.
\]
On the other hand, from Theorem \ref{themJ}
\[\beta_q=\dim C_q(J,\infty)=0, q=0,1,2\cdots.
\]
But this contradicts the Morse inequality (\ref{morse}), and hence $J$ must admit a nontrivial critical point.
i.e., there exists a nontrivial solution $u$ of
\begin{equation*}
\begin{cases}
-\Delta u+(V+m)u=f^*_m(x,u)=(a+m)u^-+g(x,u^-),& x\in \mathbb {R}^N,\\
u(x)\to 0,& as |x|\to +\infty,
\end{cases}
\end{equation*}
notice $f^*_m(x,u)\leq0$,
according to the weak maximum principle for $\mathbb{R}^N$ of \cite{12}(see Proposition 8.1 \cite{12}), the solution is negative, so it is also a negative solution of

\begin{equation*}
\begin{cases}
-\Delta u+(V+m)u=f_m(x,u)=(a+m)u^-+(b+m)u^++ g(x,u),& x\in \mathbb {R}^N,\\
u(x)\to 0,& as |x|\to +\infty,
\end{cases}
\end{equation*}
i.e., there exists a negative solution of (\ref{2}). We complete the proof of Theorem \ref{them1} thereby.\qed \vskip 5pt

\textbf{Acknowledgment. }
The first author is supported by NSFC(11871066), NSFY(Y911251) and
NSFC(11471319).

\end{document}